\input amstex\documentstyle{amsppt}  
\pagewidth{12.5cm}\pageheight{19cm}\magnification\magstep1
\NoBlackBoxes
\topmatter
\title On the characters of unipotent representations of a semisimple $p$-adic group\endtitle
\author Ju-Lee Kim and George Lusztig\endauthor
\rightheadtext{Characters of unipotent representations of a $p$-adic group}
\address{Department of Mathematics, M.I.T., Cambridge, MA 02139}\endaddress
\email{julee\@ math.mit.edu, gyuri\@ math.mit.edu}\endemail
\abstract
{Let $G$ be a semisimple almost simple algebraic group defined and split over a nonarchimedean local field 
$K$ and let $V$ be a unipotent representation of $G(K)$ (for example, an Iwahori-spherical representation). 
We calculate the character of $V$ at compact very regular elements of $G(K)$. }
\endabstract
\thanks{Both authors are supported in part by the National Science Foundation}\endthanks
\endtopmatter   
\document
\define\pos{\text{\rm pos}}
\define\Irr{\text{\rm Irr}}
\define\sneq{\subsetneqq}

\define\ucp{\un{\cp}}

\define\uK{\un K}

\define\uQ{\un Q}
\define\uB{\un B}

\define\mpb{\medpagebreak}

\define\bQ{\bar Q}
\define\bT{\bar T}

\define\si{\sim}
\define\wt{\widetilde}
\define\sqc{\sqcup}

\define\bZ{\bar Z}

\define\lb{\linebreak}

\define\op{\oplus}

\define\part{\partial}

\define\n{\notin}
\define\iy{\infty}
\define\m{\mapsto}
\define\do{\dots}

\define\lra{\leftrightarrow}

\define\sub{\subset}    

\define\T{\times}
\define\ti{\tilde}
\define\nl{\newline}
\redefine\i{^{-1}}

\define\un{\underline}

\define\ot{\otimes}

\define\Ad{\text{\rm Ad}}
\define\Hom{\text{\rm Hom}}

\define\tr{\text{\rm tr}}

\define\supp{\text{\rm supp}}

\define\a{\alpha}
\redefine\b{\beta}

\define\g{\gamma}

\redefine\o{\omega}
\define\p{\pi}
\define\ph{\phi}
\define\ps{\psi}
\define\r{\rho}
\define\s{\sigma}
\redefine\t{\tau}

\define\k{\kappa}
\redefine\l{\lambda}

\redefine\G{\Gamma}

\define\Om{\Omega}

\define\Th{\Theta}

\define\boc{\bold c}

\define\CC{\bold C}

\define\EE{\bold E}

\define\II{\bold I}

\define\NN{\bold N}

\define\ZZ{\bold Z}

\define\cc{\Cal C}

\define\ce{\Cal E}

\define\cg{\Cal G}
\define\ch{\Cal H}

\define\cl{\Cal L}
\define\cm{\Cal M}

\define\co{\Cal O}
\define\cp{\Cal P}

\define\cu{\Cal U}
\define\cv{\Cal V}
\define\cw{\Cal W}

\define\cx{\Cal X}

\define\fc{\frak c}

\define\fm{\frak m}

\define\fp{\frak p}

\define\ft{\frak t}

\define\fA{\frak A}

\define\fS{\frak S}

\define\tio{\ti\o}

\define\sha{\sharp}

\define\uP{\un P}

\define\bg{\bar g}

\define\bP{\bar P}

\define\AK{AK}
\define\BK{BK}
\define\BT{BT}
\define\CA{C1}
\define\CAA{C2}
\define\DL{DL}
\define\vD{vD}
\define\GE{Ge}
\define\HC{H1}
\define\HCC{H2}
\define\HCCC{H3}
\define\IM{IM}
\define\KL{KL}
\define\KmL{KmL}
\define\REP{L1}
\define\SQ{L2}
\define\ORA{L3}
\define\CLA{L4}
\define\MS{MS}
\define\MP{MP}
\head 1. Introduction\endhead
\subhead 1.1\endsubhead
Let $K$ be a nonarchimedean local field. Let $\co$ and $\fp$ denote the ring of integers of $K$ and the 
maximal ideal of $\co$ respectively. Let $G$ be a semisimple almost simple algebraic group defined and split
over $K$ with a given $\co$-structure compatible with the $K$-structure. 

If $V$ is an admissible representation of $G(K)$ of finite length, we denote by $\ph_V$ the character of $V$
in the sense of Harish-Chandra, viewed as a $\CC$-valued function on the set $G(K)_{rs}:=G_{rs}\cap G(K)$.
(Here $G_{rs}$ is the set of regular semisimple elements of $G$ and $\CC$ is the field of complex numbers.)
When $V$ is supercuspidal, $\ph_V$ can be sometimes studied via Harish-Chandra\rq{}s integral formula 
(\cite{\HC}). When $V$ is a fully parabolically induced representation, it is possible to compute its 
character in terms of its inducing data (see \cite{\vD}). However, when $V$ is a subquotient representation,
not much is known about the explicit values of $\ph_V$. In \cite{\KmL}, we studied the character of the 
Steinberg representation in an attempt to understand characters of subquotient representations.  For the 
Steinberg representation, van Dijk\rq{}s formula was enough to compute character values at very regular 
elements (see 1.2) thanks to the expression of the Steinberg character as a virtual character 
(\cite{\CA, \HCC}). However, in general, one needs new ideas to make a computation.

In this paper we use the local constancy of characters (\cite{\HCCC}, in more precise form given by 
\cite{AK}, \cite{\MS}) and Hecke algebras to study the characters of unipotent representations. (The notion 
of unipotent representation of $G(K)$ is defined in \cite{\SQ} (see also \cite{\CLA}) assuming that $G$ is 
adjoint; the same definition can be given without this assumption.) More precisely, we study the restriction
of the function $\ph_V$ to a certain subset $G(K)_{cvr}$ of $G(K)_{rs}$, that is to the set of compact very 
regular elements in $G(K)$ (see 1.2), in the case where $V$ is an irreducible unipotent representation of 
$G(K)$. We show that $\ph_V|_{G(K)_{cvr}}$ takes integer values which, at least if $V$ is tempered and $G$ 
is adjoint, are explicitly computable. (See Theorem 4.6). We expect that the values of $\ph_V$ on compact 
elements of $G(K)_{rs}$ which are not in $G(K)_{cvr}$ are much more complicated. Our results apply in 
particular in the case where $V$ is an irreducible admissible representation of $G(K)$ with nonzero vectors 
fixed by an Iwahori subgroup; they provide a partial generalization of the results in \cite{\KmL}.
 
In the case $g\in G(K)_{vr}$ (see 1.2) is contained in a $K$-split torus and $V$ is Iwahori spherical, 
$\ph_V(g)$ can be expressed in terms of the character of the correponding Iwahori Hecke algebra module (see 
Theorem 4.3 in \cite{\KmL}.) In general, when $g\in G(K)_{vr}$, the problem of computing the character of a 
unipotent representation $V$ of $G(K)$ at $g$ can in principle be reduced to the special case where $g$ is 
compact very regular. Indeed, the results of \cite{\CAA} reduce this problem to the computation of the 
character of a Jacquet module of $V$ (which is an admissible representation of finite length of a possibly 
smaller group with composition factors being unipotent representations computable from \cite{\BK, Cor. 8.4})
on a compact very regular element, for which the results in our paper (see Section 4) are applicable. 

\subhead 1.2\endsubhead
{\it Notation.} Let $\uK$ be a maximal unramified field extension of $K$. Let $\un\co$ and $\un\fp$ denote 
the ring of integers  of $\uK$ and the maximal ideal of $\un\co$ respectively. Let $\uK^*=\uK-\{0\}$.

Let $g\in G_{rs}\cap G(\uK)$. Let $T'=T'_g$ be the maximal torus of $G$ that contains $g$. We say that $g$ 
is very regular (resp. compact very regular) if $T'$ is split over $\uK$ and for any root $\a$ with respect 
to $T'$ viewed as a homomorphism $T'(\uK)@>>>\uK^*$ we have 
$$\a(g)\n(1+\un\fp)\ \text{ (resp.}\ \a(g)\in\un\co-\un\fp,\ \a(g)\n(1+\un\fp)).$$
Let $G(\uK)_{vr}$ (resp. $G(\uK)_{cvr}$) be the set of elements in $G(\uK)$ which are very regular (resp. 
compact very regular). Note that $G(\uK)_{cvr}\sub G(\uK)_{vr}\sub G_{rs}$. We write 
$G(K)_{vr}=G(\uK)_{vr}\cap G(K)$, $G(K)_{cvr}=G(\uK)_{cvr}\cap G(K)$.

We write $\co/\fp=F_q$, a finite field with $q$ elements, of characteristic $p$. 
Let $\uK_*$ be the set of all $x\in\uK$ such that $x^n=1$ for some $n\ge1$ prime to $p$. 

Let $\cp$ (resp. $\ucp$) be the set of parahoric subgroups (see \cite{\IM}, \cite{\BT}) of $G(K)$ (resp. of 
$G(\uK)$). Let $G(\uK)_{der}$ (resp. $G(K)_{der}$) be the derived subgroup of $G(\uK)$ (resp. $G(K)$). Let 
$G(\uK)'$ be the subgroup of $G(\uK)$ generated by $G(\uK)_{der}$ and by an Iwahori subgroup \cite{\IM,\BT}) 
of $G(\uK)$. Let $G(K)'=G(\uK)'\cap G(K)$. 

\head 2. Compact very regular elements\endhead
\subhead 2.1\endsubhead
This section contains a number of definitions and lemmas which will used in the proof of the character
formula in \S4. (The definitions in 2.3 are an exception to this; they are only used in 2.5(c) and \S5.)

Let $F:\uK@>>>\uK$ be the Frobenius automorphism, that is the field automorphism whose restriction to 
$\uK_*$ is the map $x\m x^q$; then $K$ is the fixed point set of $F$. Now $F$ induces a group isomorphism 
$G(\uK)@>>>G(\uK)$ (denoted again by $F$) whose fixed point set is $G(K)$. Note that $F$ maps any parahoric 
subgroup of $G(\uK)$ onto a parahoric subgroup of $G(\uK)$. The map $\uP\m\uP\cap G(K)$ is a $1-1$ 
correspondence between the set $\{\uP\in\ucp;F(\uP)=\uP\}$ and $\cp$.

Let $\ucp_{\emptyset}$ be the set of $\uP\in\ucp$ that are Iwahori subgroups. The minimal elements (under 
inclusion) of $\ucp-\ucp_\emptyset$ fall into finitely many $G(\uK)'$-orbits $(\ucp_i)_{i\in\II}$ (under 
conjugation); here $\II$ is a finite indexing set. If $J\sneq\II$ let $\ucp_J$ be the set of $\uP\in\ucp$ 
such that the following holds: for $i\in\II$, $\uP$ contains a subgroup in $\ucp_i$ if and only if $i\in J$.
The sets $\ucp_J$ with $J\sneq\II$ are exactly the $G(\uK)'$-orbits on the set of parahoric subgroups of 
$G(\uK)$. This agrees with the earlier notation $\ucp_\emptyset$; for $i\in\II$ we have 
$\ucp_{\{i\}}=\ucp_i$.

Let $T'$ be a maximal torus of $G$ which is defined and split over $\uK$. We define two subgroups 
$T'_*,T'_!$ of $T'(\uK)$ as follows. Let $Y'$ be the group of cocharacters of $T'$. We can identify 
$\uK\ot Y'=T'(\uK)$ in an obvious way. Let $T'_*$ (resp. $T'_!$) be the subgroup of $T'(\uK)$ corresponding 
to $\uK_*\ot Y$ (resp. $(1+\un\fp)\ot Y$) under this isomorphism; thus $T'_*$ is the set of all 
$t\in T'(\uK)$ such that $t^n=1$ for some $n\ge1$ prime to $p$. Let $T'_{*!}=T'_*T'_!$, a subgroup of 
$T'(\uK)$. 

We fix a maximal torus $T$ of $G$ defined and split over $K$ and an Iwahori subgroup $\uB$ 
of $G(\uK)$ such $F(\uB)=\uB$ and $T_{*!}\sub\uB$. 

Let $X$ be the group of characters of $T$. Let $R\sub X$ be the set of roots with respect to $T$. For each 
$\a\in R$ let $U_\a$ be the corresponding root subgroup of $G$.

Let $W$ be an indexing set for the set of $(\uB,\uB)$-double cosets 
in $G(\uK)$. We denote by $O_w$ the double coset corresponding to $w\in W$. Let 
$W'=\{w\in W;O_w\sub G(\uK)'\}$ and let $\Om=\{w\in W;O_w\sub N(\uB)\}$; here $N(\uB)$ is the normalizer of 
$\uB$ in $G(\uK)$. If $i\in\II$ let $\uP_i$ be the unique element of $\ucp_i$ containing $\uB$. Then 
$\uP_i-\uB=O_w$ for a well defined element $w\in W'$; we set $w=s_i$. 
For $J\sneq\II$ let $\uP_J$ be the unique subgroup in $\ucp_J$ that contains $\uB$.
Let $W_J$ be the subgroup of $W'$ generated by $\{s_i;i\in J\}$.

Two $(\uB,\uB)$-double cosets 
$O_w,O_{w'}$ are said to be composable if there exists $w''\in W$ such that multiplication defines a 
bijection $O_w\T_{\uB}O_{w'}@>>>O_{w''}$; here $\uB$ acts on $O_w\T O_{w'}$ by $b:(g,g')\m(gb\i,bg')$. 
We then set $w\cdot w'=w''$. There is a unique group structure on $W$ such that the product of $w,w'$ is
$w\cdot w'$ whenever $O_w,O_{w'}$ are composable and $s_i^2=1$ for any $i\in\II$.
(The unit element is the $(\uB,\uB)$-double coset $\uB$.) Then $W'$ is the subgroup of $W$ generated by
$\{s_i;i\in\II\}$. This is a normal subgroup of $W$ and $\Om$ is an abelian subgroup of $W$ which maps 
isomorphically onto $W/W'$. The group $W'$ together with $\{s_i;i\in\II\}$ is a Coxeter group (in fact, an 
affine Weyl group) with length function $l:W'@>>>\NN$. We extend this to a function
$l:W@>>>\NN$ by $l(w_1w_2)=l(w_2w_1)=l(w_1)$ for $w_1\in W',w_2\in\Om$. Then $O_w,O_{w'}$ are composable
precisely when $l(ww')=l(w)+l(w')$.

Now any $\uP\in\ucp$ has a prounipotent radical $U_{\uP}$ (a normal subgroup) such that $\uP/U_{\uP}$ is 
naturally a connected reductive group over $\un\co/\un\fp$ (an algebraic closure of $F_q$).

The orbits of the diagonal conjugation action of $G(\uK)'$ on $\ucp_\emptyset\T\ucp_\emptyset$ are indexed 
by the elements of $W'$: to $w\in W'$ corresponds the $G(\uK)'$-orbit of $(\uB,g\uB g\i)$ where $g$ is 
some/any element of $O_w$. We write $\pos(\uB',\uB'')=w$ when the $G(\uK)'$-orbit of 
$(\uB',\uB'')\in\ucp_\emptyset\T\ucp_\emptyset$ corresponds to $w\in W'$.

Let $g\in G(\uK)_{cvr}$. Let $T'=T'_g$; then $T'$ is defined and split over $\uK$. We have $g\in T'_{*!}$ 
hence $g=g_*g_!$ where $g_*\in T'_*$, $g_!\in T'_!$. From the definitions we see that $g_*\in G(\uK)_{cvr}$.
Let $\ucp_\emptyset^g$ (resp. $\ucp_\emptyset^{g_*}$) be the set of $\uB'\in\ucp_\emptyset$ such that 
$g\in\uB'$ (resp. $g_*\in\uB'$). 
Let $\ucp_\emptyset^{T'_*}=\{\uB'\in\ucp_\emptyset;T'_*\sub\uB'\}$,
$\ucp_\emptyset^{T'_{*!}}=\{\uB'\in\ucp_\emptyset;T'_{*!}\sub\uB'\}$.

\proclaim{Lemma 2.2} (a) With the notation above we have 
$\ucp_\emptyset^g=\ucp_\emptyset^{g_*}=\ucp_\emptyset^{T'_*}=\ucp_\emptyset^{T'_{*!}}$.

(b) Given $w\in W'$ and $\uB'\in\ucp_\emptyset^{T'_*}=\ucp_\emptyset^g$ there is a unique subgroup 
$\uB''\in\ucp_\emptyset^{T'_*}=\ucp_\emptyset^g$ such that $\pos(\uB',\uB'')=w$. We set $\uB''=w\circ\uB'$. 
Then $w:\uB'\m w\circ\uB'$ defines a free transitive action of $W'$ on the set 
$\ucp_\emptyset^{T'_*}=\ucp_\emptyset^g$.
\endproclaim
We can assume that $T'=T$. 

For any $w\in W'$ and any $\uB_0\in\cp_\emptyset^{g_*}$ let 
$\cp_\emptyset(w,\uB_0)=\{\uB'\in\cp_\emptyset;\pos(\uB_0,\uB')=w\}$. We first show:

(c) For any $w\in W'$, $\ucp_\emptyset^{g_*}\cap\ucp_\emptyset(w,\uB_0)$ consists of exactly one element.
\nl
We argue by induction on $l(w)$. When $w=1$ the result is obvious; we have
$\ucp_\emptyset^{g_*}\cap\ucp_\emptyset(1,\uB_0)=\{\uB_0\}$.
Assume that (c) is known when $l(w)=1$. Let $w\in W'-\{1\}$. We can find
$w_1,s\in W'$ such that $w=w_1s$, $l(w)=l(w_1)+l(s)$, $l(s)=1$. If 
$\uB_1,\uB_2\in\ucp_\emptyset^{g_*}\cap\ucp_\emptyset(w,\uB_0)$ then there are unique
$\uB'_1,\uB'_2\in\ucp_\emptyset$ such that $\pos(\uB_0,\uB'_1)=w_1$, $\pos(\uB'_1,\uB_1)=s$, 
$\pos(\uB_0,\uB'_2)=w_1$, $\pos(\uB'_2,\uB_2)=s$. Since $g_*\uB_0g_*\i=\uB_0$, $g_*\uB_1g_*\i=\uB_1$, we 
must have $g_*\uB'_1g_*\i=\uB'_1$ (by uniqueness of $\uB'_1$) and similarly $g_*\uB'_2g_*\i=\uB'_2$. Thus 
$\uB'_1,\uB'_2\in\ucp_\emptyset^{g_*}\cap\ucp_\emptyset(w_1,\uB_0)$. By the induction hypothesis we have 
$\uB'_1=\uB'_2$ (which we denote by $\b$). We have 
$\uB_1,\uB_2\in\ucp_\emptyset^{s_*}\cap\ucp_\emptyset(s,\b)$. Since (c) is assumed to hold for $s$ instead 
of $w$ we see that $\uB_1=\uB_2$. Thus $\ucp_\emptyset^{g_*}\cap\ucp_\emptyset(w,\uB_0)$ consists of at most
one element. By the induction hypothesis we can find 
$\b'\in\ucp_\emptyset^{g_*}\cap\ucp_\emptyset(w_1,\uB_0)$ and by our assumption we can find
$\uB'\in\ucp_\emptyset^{g_*}\cap\ucp_\emptyset(s,\b')$. Then we have
$\uB'\in\ucp_\emptyset^{g_*}\cap\ucp_\emptyset(w,\uB_0)$. We see that (c) holds for $w$. It remains to
verify (c) in the case where $w=s$ is of length $1$.
Let $\b\in\cp_\emptyset(s,\uB_0)$. We can assume that $\b$ contains $T_*$ hence also $g_*$.
(We can take $\b$ to be a conjugate of $B_0$ under an element of $G(\uK)'$ normalizing $T_*$.)
We can find an injective group homomorphism $f:\co@>>>G(\uK)'$ such that 
$t\m f(t)\b f(t)\i$ is a bijection $\uK_*\cup\{0\}@>\si>>\cp_\emptyset(s,\uB_0)$ and
$gf(t)g\i=f(x_0t)$ for all $t\in\uK_*\cup\{0\}$ for some $x_0\in\uK_*-\{1\}$.
Then $\ucp_\emptyset^{g_*}\cap\ucp_\emptyset(s,\uB_0)$ is in bijection with
$Z:=\{t\in\uK_*\cup\{0\};gf(t)\b f(t)\i g\i=f(t)\b f(t)\i\}=\{t\in\uK_*\cup\{0\};f(t)\i gf(t)\in\b\}$.
For $t\in\uK_*\cup\{0\}$ we have $f(t)\i gf(t)\in\b$ if and only if $f(t)\i gf(t)g\i=f(-t+x_0t)\in\b$ (since
$g\in\b$) that is if and only if $f(-t+x_0t)=1$ (we use that $f(\co)\cap\b=\{1\}$ since $f$ is injective) 
that is if and only if $t(x_0-1)=0$ which is the same as $t=0$ (since $x_0\ne0$). We see that $Z=\{0\}$.
Hence $\ucp_\emptyset^{g_*}\cap\ucp_\emptyset(s,\uB_0)$ consists of exactly one element.
This completes the proof of (c).

We shall denote the unique element in (c) by $w\circ\uB_0$. Let $\fA=\{w\circ\uB;w\in W'\}$. From (c) we 
have 

$\ucp_\emptyset^{g_*}=\fA$.
\nl
Let $\uB'\in\ucp_\emptyset^{g_*}$ and let $t\in T_{*!}$. Since $tg_*=g_*t$ we have 
$t\uB't\i\in\ucp_\emptyset^{g_*}$.
Thus if $\uB'\in\fA$ then $t\uB't\i\in\fA$. We have $\uB'\in\cp_\emptyset(w,\uB)$ for some $w\in W'$
and $t\uB't\i\in\cp_\emptyset(w,\uB)$ (since $t\in\uB$). Thus both
$\uB',t\uB't\i$ belong to $\ucp_\emptyset^{g_*}\cap\ucp_\emptyset(w,\uB)$; hence, by (c) we have
$\uB'=t\uB't\i$ so that $t\in\uB'$. We see that $\uB'\in\ucp_\emptyset^{T_{*!}}$. Thus
$\fA\sub\ucp_\emptyset^{T_{*!}}$.

We can find $n\ge1$ such that $g_*^{p^n}=g_*$. We have 
$\lim_{m\to\iy}g_!^{p^{nm}}=1$ hence $\lim_{m\to\iy}g^{p^{nm}}=g_*$. If $B''\in\ucp_\emptyset^g$, then 
$g^{p^{nm}}B''g^{-p^{nm}}=B''$ for any $m\ge1$. Taking $m\to\iy$ we deduce $g_*B''g_*\i=B''$ hence 
$B''\in\ucp_\emptyset^{g_*}$. Thus $\ucp_\emptyset^g\sub\ucp_\emptyset^{g_*}$. The inclusion
$\ucp_\emptyset^{T_{*!}}\sub\ucp_\emptyset^g$ is obvious. Thus we have 
$$\ucp_\emptyset^{T_{*!}}\sub\ucp_\emptyset^g\sub\ucp_\emptyset^{g_*}=\fA\sub\ucp_\emptyset^{T_{*!}},$$
hence
$$\ucp_\emptyset^{T_{*!}}=\ucp_\emptyset^g=\ucp_\emptyset^{g_*}=\fA.$$
We have clearly 
$$\ucp_\emptyset^{T_{*!}}\sub\ucp_\emptyset^{T_*}\sub\ucp_\emptyset^{g_*}=\ucp_\emptyset^{T_{*!}}$$
hence $\ucp_\emptyset^{T_{*!}}=\ucp_\emptyset^{T_*}$.
This completes the proof of (a).

Now (b) follows immediately from (a) and (c). (Note that $\ucp_\emptyset^{T'_*}$ is nonempty; it 
contains $\uB$.) The lemma is proved.

\subhead 2.3\endsubhead
Let $\wt{G(\uK)_{cvr}}$ be the set of all pairs $(g,\uB')$ where $g\in G(\uK)_{cvr}$ and 
$\uB'\in\ucp_\emptyset^g$. By 2.2(b), $w:(g,\uB')\m(g,w\circ\uB')$ is a free action of $W'$ on 
$\wt{G(\uK)_{cvr}}$ whose orbits are exactly the fibres of the first projection 
$pr_1:\wt{G(\uK)_{cvr}}@>>>G(\uK)_{cvr}$. 

Now let $A$ be a (finite dimensional) representation of $W'$. For any $g\in G(\uK)_{cvr}$ we define a vector
space $A_g$ as the set of the $W'$-orbits on the set $\ucp_\emptyset^g\T A$ for the free $W'$-action
$w:(\uB',a)\m(w\circ\uB',wa)$. We can think of the union $\sqc_{g\in G(\uK)_{cvr}}A_g$ and its obvious
projection to $G(\uK)_{cvr}$ as something like a "local system" 
$\cl_A$ on $G(\uK)_{cvr}$ associate to $A$ and the principal $W'$-covering $pr_1$ above.

We now define an isomorphism of local systems $\ps:F^*\cl_A@>\si>>\cl_A$ as the collection 
of isomorphisms $\ps_g:A_{F(g)}@>>>A_g$ (with $g\in G(\uK)_{cvr}$) where $\ps_g$ is induced by

$\ucp_\emptyset^{F(g)}\T A@>>>\ucp_\emptyset^g\T A$, $(\uB',a)\m(F\i(\uB'),a)$
\nl
by passage to $W'$-orbits.

We define $\t_A:G(K)_{cvr}@>>>\CC$ by $\t_A(g)=\tr(\ps_g:A_g@>>>A_g)$. 

\proclaim{Lemma 2.4} Let $g\in G(\uK)_{cvr}$. Let $T'=T'_g$ and $g_*$ be as in 2.1. Let $\uP\in\ucp$.

(a) The following four conditions are equivalent: (i) $g\in\uP$; (ii) $g_*\in\uP$; (iii) $T'_*\sub\uP$; 
(iv) $T'_{*!}\sub\uP$.

(b) If $g\in\uP$ and $\bg$ is the image of $g$ in $\uP/U_{\uP}$ then 
$\bg$ is regular semisimple in $\uP/U_{\uP}$.
\endproclaim

We prove (a). Assume that $g\in\uP$. We can find a Borel subgroup $\b$ of $\uP/U_{\uP}$ which contains the 
image of $\bg$ of $g$ in $\uP/U_{\uP}$. The inverse image of $\b$ under $\uP@>>>\uP/U_{\uP}$ is an Iwahori 
subgroup $\uB'$ of $G(\uK)$ which contains $g$ and is contained in $\uP$. By 2.2(a), we have 
$T'_{*!}\sub\uB'$ hence $T'_{*!}\sub\uP$. Thus (i) implies (iv). An entirely similar proof shows that (ii) 
implies (iv). It is obvious that (iv) implies (iii), that (iii) implies (ii) and that (iv) implies (i). This
 proves (a).

We prove (b). The image of $T'_*$ under $\uP@>>>\uP/U_{\uP}$ is a maximal torus $\bT'$ of $\uP/U_{\uP}$ that
contains $\bg$. For any root $\a$ of $\uP/U_{\uP}$, viewed as a character $\bT'@>>>\un\co/\un\fp-\{0\}$, the
value $\a(\bg)$ is equal to $\a'(g_*)\in\uK_*$ for some root $\a'$ of $G$ viewed as a character of $T'$. (We 
identify $\uK_*=\un\co/\un\fp-\{0\}$ in an obvious way.) Since $\a'(g_*)\ne1$ we must have $\a(\bg)\ne1$. 
This proves (b).

\proclaim{Lemma 2.5} Let $T'$ be a maximal torus of $G$ which is defined over $K$ and is $\uK$-split. 

(a) There is a unique $W'$-conjugacy class $C_{T'}$ in $W'$ such that the following holds: for any
$\uB'\in\ucp_\emptyset^{T'_*}$ we have $\pos(\uB',F(\uB'))=w$ for some $w\in C_{T'}$.

(b) Any element of $C_{T'}$ has finite order in $W'$.

(c) Let $A,\t_A$ be as in 2.3 and let $g\in G(K)_{cvr}\cap T'$. We have $\t_A(g)=\tr(w\i,A)=\tr(w,A)$
where $w\in C_{T'}$. Moreover, we have $\t_A(g)\in\ZZ$. 
\endproclaim
We prove (a). It is enough to show that if $\uB'\in\ucp_\emptyset^{T'_*}$, $\uB''\in\ucp_\emptyset^{T'_*}$ 
and $w',w''$ in $W'$ are given by $\pos(\uB',F(\uB'))=w'$, $\pos(\uB'',F(\uB''))=w''$ then $w''=yw'y\i$ for 
some $y\in W'$. We have $F(\uB')\in\ucp_\emptyset^{T'_*}$, $F(\uB'')\in\ucp_\emptyset^{T'_*}$. From the 
definitions we have $F(\uB')=w'\circ\uB'$, $F(\uB'')=w''\circ\uB''$. We have $\uB''=y\circ\uB'$ for a 
unique $y\in W'$. It follows that $F(\uB'')=y\circ F(\uB')$. (Note that $F$ induces the identity map on 
$W'$.) Thus $w''\circ(y\circ\uB')=y\circ(w'\circ\uB')$ that is $(w''y)\circ\uB'=(yw')\circ\uB'$. By the 
freeness of the $W'$-action on $\ucp_\emptyset^{T'_*}$ we then have $w''y=yw'$. This proves (a).

We prove (b). Let $\uB'\in\ucp_\emptyset^{T'_*}$. We set $\pos(\uB',F(\uB'))=w$ so that $w\in C_{T'}$. Now 
$T'$ becomes split after a finite unramified extension of $K$. Hence there exists $s\ge1$ such that 
$F^s(\uB')=\uB'$. We have $F(\uB')=w\circ\uB'$, $F^2(\uB')=w\circ(F(\uB'))=w^2\circ\uB'$, $\do$, 
$F^s(\uB')=w^s\circ\uB'$. Thus $w^s\circ\uB'=\uB'$ so that $w^s=1$. This proves (b).

The first assertion of (c) follows from definitions. To prove the second assertion of (c) we can assume that
$w\in W_H$ for some $H\sneq\II$ so that $\tr(w,A)=\tr(w,A|_{W_H})$. We then use the rationality of the
irrducible representations of $W_H$. The lemma is proved.

\subhead 2.6\endsubhead
In \cite{\GE, 3.3}, G\'erardin defines a map 
$$\align&\mu:\{G(K)_{der}-\text{conjugacy classes of }\uK-\text{split maximal tori in }G
\text{ defined over }K\}
@>>>\\&\{W'-\text{conjugacy classes of elements of finite order in }W'\}\endalign$$
and shows that

(a) $\mu$ is a bijection.
\nl
(Actually in \cite{\GE} it is assumed that $G$ is simply connected so that $G(K)_{der}=G(K)$ but the general 
case can be deduced from this.) One can show that $\mu$ is induced with the map induced by $T'\m C_{T'}$ in 
2.5.

We now introduce some notation which will be used in the following lemma. 
Let $T''$ be a maximal torus of $G$ which is defined over $K$ and is split over $\uK$. Using the obvious 
identification $\uK_*=\un\co/\un\fp-\{0\}$ we can view $T''_*$ as a torus over $\un\co/\un\fp$ (identified 
with its group of $\un\co/\un\fp$ points) with a natural $F_q$-structure (induced by $F$). Hence $T''_*$
has a maximal $F_q$-split subtorus $T''_{*s}$ (it is the identity component of $\{t\in T''_*;t^q=t\}$). Now 
let $g\in G(K)_{cvr}$. Let $T'=T'_g$ and let $C=C_{T'}$. Using 2.5(b) we see that we can find $H\sneq\II$ 
and $w\in W_H$ such that $w\in C$ and $w$ is elliptic in $W_H$. 

\proclaim{Lemma 2.7} (a) There exists $\uP\in\ucp_H$ such that $F(\uP)=\uP$ and $T'_*\sub\uP$ (hence 
$g\in\uP$, see 2.4(a)).

(b) There exists a $G(K)_{der}$-conjugate $g'$ of $g$ such that $g'\in\uP_H$ and the maximal 
$(T'_{g'})_{*s}\sub T_*$.
\endproclaim
Let $\uB''\in\ucp_\emptyset^{T'_*}$. Then $F(\uB'')=w'\circ\uB''$ for a unique $w'\in C$. We have 
$w'=ywy\i$ for some $y\in W'$. Setting $\uB'=y\i\circ\uB''$ we have $\uB'\in\ucp_\emptyset^{T'_*}$, 
$F(\uB')=w\circ\uB'$  hence $\pos(\uB',F(\uB'))=w$. Define $\uP\in\cp_H$ by $\uB'\sub\uP$. Then we have 
automatically $F(\uB')\sub\uP$ (since $\pos(\uB',F(\uB'))\in W_H$). We have also $F(\uB')\sub F(\uP)$. Since
$F(\uB')$ is contained in a unique parahoric subgroup in $\ucp_H$ we must have $F(\uP)=\uP$. We have 
$T'_*\sub\uP$ and (a) is proved.

We prove (b). Let $\uP$ be as in (a). The image of $T'_{*s}$ under $pr:\uP@>>>\uP/U_{\uP}$ is a torus of 
$\uP/U_{\uP}$ defined and split over $F_q$ and we can find a Borel subgroup $\b$ of $\uP/U_{\uP}$ defined 
over $F_q$ that contains this torus. Now $pr\i(\b)$ is an Iwahori subgroup $\uB'$ of $G(\uK)$ such that 
$\uB'\sub\uP$, $F(\uB')=\uB'$ and $T'_{*s}\sub\uB'$. We can find a maximal torus $T''$ of $G$ defined and 
split over $K$ such that $T''_*\sub\uB'$. Now $T''_*$ viewed as a torus over $\un\co/\un\fp$ is an 
$F_q$-split torus in $\uB'$. We can find $x\in\uB'\cap G(K)$ such that $T'_{*s}\sub xT''_*x\i$. Replacing 
$T''$ by $xT''x\i$ we can assume that $T'_{*s}\sub T''_*\sub\uB'$. We can find $z\in G(K)_{der}$ such that 
$zT''z\i=T,z\uB' z\i=\uB$. Let $g'=zgz\i,\uP'=z\uP z\i$. We have $\uB\sub\uP'$, $\uP'\in\ucp_H$ hence 
$\uP'=\uP_H$. We have $T'_{g'}=zT'z\i$, $(zT'z\i)_{*s}\sub T_*$, $g'\in\uP'$. This proves (b).

\proclaim{Lemma 2.8} Assume that $g\in G(K)_{cvr},g'\in G(K)_{cvr}$ are contained in the same maximal torus 
$T'$ of $G$. If $\uP$ is a minimal $F$-stable parahoric subgroup of $G(\uK)$ containing $g$ then
$\uP$ is a minimal $F$-stable parahoric subgroup of $G(\uK)$ containing $g'$.
\endproclaim
This follows from the following statement.

If $\uP\in\ucp$, then the following two conditions are equivalent: (i) $g\in\uP$ and (ii) $g'\in\uP$. 
\nl
Indeed, by 2.4(a) both conditions are equivalent to the condition that $T'_*\sub\uP$.

\head 3. Recollections on affine Hecke algebras\endhead
\subhead 3.1\endsubhead
Our proof of the character formula in \S4 will involve detailed information on the affine Hecke algebras 
which appear in the study of unipotent representations. In preparation for that proof, we now review some 
definitions and results of \cite{\CLA}. (In \cite{\CLA} , $G$ is assumed to be adjoint but the results of 
\cite{\CLA} that we use in this section extend with the same proof to the general case.)
 
We set $B:=\uB\cap G(K)$; this is an Iwahori subgroup of $G(K)$.

For any $J\sneq\II$ let $\cp_J$ be the set of $P\in\cp$ of the form $\uP\cap G(K)$ where $\uP\in\ucp_J$ is 
$F$-stable. The sets $\cp_J$ are exactly the orbits of $G(K)'$ on $\cp$ (under conjugation). Let $P_J$ be 
the unique subgroup in $\cp_J$ that contains $B$. 

Now any $P\in\cp$ has a prounipotent radical $U_P$ (a normal subgroup) such that the quotient $\bP=P/U_P$ is
naturally the group of $F_q$-points of a connected reductive group defined and split over $F_q$.

\subhead 3.2\endsubhead
We now fix $J\sneq\II$ and a unipotent cuspidal (irreducible) representation (over $\CC$) of the finite 
group $\bP$ where $P=P_J$; let $\EE$ be the vector space (over $\CC$) of this representation. We regard 
$\EE$ as a $P$-module via the surjective homomorphism $P@>>>\bP$.

Let $\cx$ be the $\CC$-vector space consisting of all functions $f:G(K)@>>>\EE$ such that $f(gh)=h\i f(g)$ 
for all $g\in G(K),h\in P$ and $\supp(f)=\{g\in G(K);f(g)\ne0\}$ is contained in a union of finitely many 
$P$-cosets in $G(K)/P$. For $g'\in G(K),f\in\cx$ we define $g'f:G(K)@>>>\EE$ by $(g'f)(g)=f(g'{}\i g)$. Then
$g'f\in\cx$. This defines a representation of $G(K)$ on the vector space $\cx$.
We have a direct sum decomposition $\cx=\op_{gP}\cx_{gP}$ where $gP$ runs over the cosets $G(K)/P$
and $\cx_{gP}=\{f\in\cx;\supp(f)\in gP\}$. Let $\ch$ be the endomorphism algebra of the $G(K)$-module $\cx$. 
For any $(P,P)$-double coset $\Th$ in $G(K)$ we denote by $\ch_\Th$ the set consisting of those $\ps\in\ch$ 
such that $\ps(\cx_{gP})\sub\op_{g'P\in G(K)/P;g'{}\i g\in\Th}\cx_{g'P}$ for some $gP\in G(K)/P$ or 
equivalently any $gP\in G(K)/P$. We say that $\Th$ (as above) is good if for some (or equivalently any) 
$g\in\Th$ we have $U_P(P\cap gPg\i)=P$. If $\Th$ is a good double coset then $\dim\ch_\Th=1$; if $\Th$ is 
not a good double coset then $\ch_\Th=0$. Hence $\ch=\op_\Th\ch_\Th$ where $\Th$ runs over the good 
$(P,P)$-double cosets in $G(K)$.

\subhead 3.3\endsubhead
Let $\Irr(G(K);P,\EE)$ be the set of irreducible admissible representations $(V,\s)$ of $G(K)$ (up to 
isomorphism) with the following property: the (necessarily finite dimensional) vector space $V^{U_P}$ of 
$U_P$-invariants, regarded as a $\bP$-module in an obvious way, contains the $\bP$-module $\EE$. An 
irreducible admissible representation of $G(K)$ is said to be unipotent if it belongs to $\Irr(G(K);P,\EE)$ 
for some $P,\EE$ as in 3.2.

We now fix again $P,\EE$ as in 3.2 and we fix $(V,\s)\in\Irr(G(K);P,\EE)$. By a version of Frobenius 
reciprocity we have an isomorphism $\Hom_{G(K)}(\cx,V)@>\si>>\Hom_P(\EE,V)$ induced by the imbedding 
$\EE@>>>\cx$ which takes $e\in\EE$ to the function $f_e:G(K)@>>>\EE$ given by $h\m h\i e$ if $h\in P$ and 
$g\m 0$ if $g\n P$. Since $\Hom_{G(K)}(\cx,V)$ is an $\ch$-module in an obvious way, we see that 
$\Hom_P(\EE,V)$ is naturally an $\ch$-module. Clearly, we have $\Hom_P(\EE,V)=\Hom_{\bP}(\EE,V^{U_P})$. 
Hence $\Hom_{\bP}(\EE,V^{U_P})$ (a finite dimensional $\CC$-vector space) becomes an $\ch$-module. 

If $\Th$ is a $(P,P)$-double coset in $G(K)$ then $\Th$ is a union $\cup_w(O_w\cap G(K))$ where $w$ runs 
over a $(W_J,W_J)$-double coset in $W$. This gives a $1-1$ correspondence between the set of 
$(P,P)$-double cosets in $G(K)$ and the set of $(W_J,W_J)$-double cosets in $W$. A $(P,P)$-double coset is 
good if and only if the corresponding $(W_J,W_J)$-double coset in $W$ is contained in $NW_J$, the normalizer
of $W_J$ in $W$. In each $(W_J,W_J)$-double coset in $W$ contained in $NW_J$ there is a unique element of 
minimal length. Hence the set $\cw$ consisting of the elements of minimal length in the various 
$(W_J,W_J)$-double cosets in $W$ contained in $NW_J$ is an indexing set for the set of good $(P,P)$-double 
cosets in $G(K)$. We denote the good $(P,P)$-double coset corresponding to an element $w\in\cw$ by $\Th_w$. 
Actually, $\cw$ is a subgroup of $NW_J$. We have $\ch=\op_{w\in\cw}\ch_{\Th_w}$ with $\dim\ch_{\Th_w}=1$ for
any $w\in\cw$. As in \cite{\CLA} we see that this decomposition gives rise to a description of $\ch$ as an 
(extended) affine Hecke algebra with explicitly known (possibly unequal) parameters.

\subhead 3.4\endsubhead
Let $H\sneq\II$ be such that $J\sub H$ and let $Q=P_H$, $\uQ=\uP_H$. Now $\b:=B/U_Q$ is the group of 
$F_q$-points of a Borel subgroup defined over $F_q$ of the reductive group $\uQ/U_{\uQ}$ defined and split 
over $F_q$ and $\p:=P/U_Q$ is the group of $F_q$-points of a parabolic subgroup of that reductive group with
$\b\sub\p$. Let $\k:Q@>>>\bQ$ be the canonical homomorphism. The subsets $\k(O_w\cap G(K))$ $(w\in W_H)$ are
exactly the $(\b,\b)$-double cosets in $\bQ$ hence $W_H$ can be identified with the Weyl group of 
$\uQ/U_{\uQ}$. For any parabolic subgroup $\p'$ of $\bQ$ we denote by $u_{\p'}$ the unipotent radical of 
$\p'$.

Let $\cx_H$ be the $\CC$-vector space consisting of all functions $f:\bQ@>>>\EE$ such that $f(gh)=h\i f(g)$ 
for all $g\in\bQ,h\in\p$. (Note that $\EE$ is naturally a $\p$-module on which $u_\p:=U_P/U_Q$ acts 
trivially.)

For $g'\in\bQ,f\in\cx_H$ we define $g'f:\bQ@>>>\EE$ by $(g'f)(g)=f(g'{}\i g)$. Then $g'f\in\cx_H$. This 
defines a representation of $\bQ$ on the vector space $\cx_H$. We have a direct sum decomposition 
$\cx_H=\op_{g\p}\cx_{H,g\p}$ where $g\p$ runs over the cosets $\bQ/\p$ and 
$\cx_{H,g\p}=\{f\in\cx_H;\supp(f)\in g\p\}$. Let $\ch_H$ be the endomorphism algebra of the $\bQ$-module 
$\cx_H$. If $f\in\cx_H$ we define $f'\in\cx$ by $f'(g)=f(\k(g))$ if $g\in Q$ and $f'(g)=0$ if $g\in G(K)-Q$.
We regard $\cx_H$ as a subspace of $\cx$ via $f\m f'$. Let $\ch'=\op_\Th\ch_\Th\sub\cx$ where $\Th$ runs 
over the good $(P,P)$-double cosets in $G(K)$ that are contained in $Q$. This is a subalgebra of $\ch$ which
can be identified with $\ch_H$ in such a way that the $\ch'$-module structure on $\cx_H$ (restriction of the
$\ch$-module structure on $\cx$) coincides with the obvious $\ch_H$-module structure on $\cx_H$. Thus 
$\ch_H$ may be identified with a subalgebra of $\ch$ (namely $\ch'$).

\subhead 3.5\endsubhead
The finite dimensional vector space $V^{U_Q}$ of $U_Q$-invariants on $V$ is naturally a $\bQ$-module (since 
$Q$ normalizes $U_Q$). This $\bQ$-module can be decomposed as $E_1\op E_2\op\do\op E_r$ where $E_1,\do,E_r$ 
are irreducible $\bQ$-modules. We show:

(a) {\it For any $i\in[1,r]$, the space of invariants $E_i^{U_P/U_Q}$ is nonzero and (when regarded as a
$P/U_P$-module in the obvious way) is a direct sum of copies of $\EE$.}
\nl

\

By a known property (see \cite{\REP, 3.25}) of unipotent representations of finite reductive groups (applied 
to the $Q/U_Q$-module $E_i$) it is enough to show that the $\p/u_\p$-module $E_i^{u_\p}$ contains a copy of 
$\EE$ (it is then nonzero and a direct sum of copies of $\EE$).

We can find a parabolic subgroup $\p'$ of $Q/U_Q$ and an irreducible cuspidal representation $\EE_0$ of 
$\p'/u_{\p'}$ such that the $\p'/u_{\p'}$-module $E_1^{u_{\p'}}$ contains $\EE_0$. Then the 
$\p'/u_{\p'}$-module $(V^{U_Q})^{u_{\p'}}$ contains $\EE_0$. Let $P'$ be the inverse image of $\p'$ under 
$Q@>>>Q/U_Q$; then $P'\in\cp$ is contained in $Q$ and the $P'/U_{P'}$-module $(V^{U_Q})^{U_{P'}}$ contains 
$\EE_0$. (We have $P'/U_{P'}=\p'/u_{\p'}$.) Replacing $P'$ by a $Q$-conjugate we can assume that
$P'=P_{J_1}$ for some $J_1\sub H$. We have $U_Q\sub U_{P'}$ hence $V^{U_{P'}}\sub V^{U_Q}$ and  
$(V^{U_Q})^{U_{P'}}=V^{U_{P'}}$. Thus the $P'/U_{P'}$-module $V^{U_{P'}}$ contains $\EE_0$. By the 
uniqueness of cuspidal support of $V$ we see that $\EE_0$ is a unipotent (cuspidal) representation and then 
using \cite{\CLA, 1.6(b)} we see that we can find $g'\in G(K)$ such that $\Ad(g')$ carries $P$ to $P'$ and 
$\EE$ to a representation isomorphic to $\EE_0$. In particular we have $\sha(J)=\sha(J_1)$. We also see that
the reductive quotients of the parabolic subgroups $\p$ and $\p'$ of $Q/U_Q$ both admit unipotent cuspidal 
representations. By the classification of unipotent cuspidal representations \cite{\ORA}, we see that these 
two parabolic subgroups are conjugate under $Q/U_Q$. Hence $P,P'$ are conjugate under $Q$ and in particular 
under $G(K)'$. Since $P=P_J,P'=P_{J_1}$ it follows that $P=P'$ and $g'$ above is such that $\Ad(g')$ carries
$P$ to $P'=P$ and $\EE$ to a representation isomorphic to $\EE_0$. By \cite{\CLA, 1.16}, $\Ad(g')$ carries 
$\EE$ to a representation isomorphic to $\EE$. Thus $\EE_0$ is isomorphic to $\EE$. We see that the 
$P/U_P$-module $E_i^{U_P/U_Q}$ contains a copy of $\EE$. This completes the proof of (a).

Replacing $Q$ by $P$ in (a) we deduce:

\smallskip

(b) {\it The (finite dimensional) $P/U_P$-module $V^{U_P}$ is $\EE$-isotypic.}

\subhead 3.6\endsubhead
We preserve the setup of 3.4.

Let $\cc_1$ be the category of (finite dimensional) representations $\cv$ of $Q/U_Q$ such that the 
$\p/u_\p$-module $\cv^{u_\p}$ of $u_\p$-invariants is a direct sum of copies of $\EE$ and it generates the 
$Q/U_Q$-module $\cv$. Note that each irreducible object in $\cc_1$ is a unipotent representation of $Q/U_Q$ 
which belongs to $\cc_1$. Let $\cc_2$ be the category of (finite dimensional) representations of the Hecke 
algebra $\ch_H$. Let $\cv\in\cc_1$. By Frobenius reciprocity we have an isomorphism
$\Hom_{Q/U_Q}(\cx_H,\cv)@>\si>>\Hom_\p(\EE,\cv)$ induced by the imbedding $\EE@>>>\cx_H$ which takes
$e\in\EE$ to the function $f_e:Q/U_Q@>>>\EE$ given by $h\m h\i e$ if $h\in\p$ and $g\m 0$ if $g\n\p$. Since 
$\Hom_{Q/U_Q}(\cx_H,\cv)$ is an $\ch_H$-module in an obvious way, we see that $\Hom_\p(\EE,\cv)$ is 
naturally an $\ch_H$-module. Clearly, we have $\Hom_\p(\EE,\cv)=\Hom_{\p/u_\p}(\EE,\cv^{u_\p})$. Hence 
$\cv_\EE:=\Hom_{\p/u_\p}(\EE,\cv^{u_\p})$ is an object of $\cc_2$. Note that $\cv\m\cv_\EE$ is an 
equivalence of categories.

For each $D\in\cc_2$ we denote by $[D]$ the object of $\cc_1$ (well defined up to isomorphism) such that 
$[D]_\EE\cong D$ in $\cc_2$.

By 3.5(a), we have $V^{U_Q}\in\cc_1$ and 
$$(V^{U_Q})_\EE=\Hom_{\p/u_\p}(\EE,V^{U_Q})=\Hom_{\p/u_\p}(\EE,(V^{U_Q})^{u_\p})=\Hom_P(\EE,V^{U_P})\in\cc_2
$$
hence 
$$V^{U_Q}\cong[\Hom_P(\EE,V^{U_P})].$$
From the definitions we see that the $\ch_H$-module $\Hom_P(\EE,V^{U_P})$ is the restriction of the
$\ch$-module $\Hom_P(\EE,V^{U_P})$ (see 3.3) to the subalgebra $\ch_H$. We see that:

(a) {\it The $Q/U_Q$-module $V^{U_Q}$ is isomorphic to a direct sum of unipotent representations
$[D_1]\op[D_2]\op\do\op[D_r]$ where $D_1,D_2,\do,D_r$ are irreducible $\ch_H$-modules such that the 
$\ch_H$-module $D_1\op D_2\op\do\op D_r$ is isomorphic to the restriction of the $\ch$-module 
$\Hom_P(\EE,V^{U_P})$ to the subalgebra $\ch_H$ of $\ch$.}

\subhead 3.7\endsubhead
Let $P=P_J$, $\EE$ be as in 3.2 and let $(V,\s)$ be as in 3.3. Let $H\sneq\II$, $Q=P_H$, $\uQ=\uP_H$ be as 
in 3.4. Let $\ce$ be the $\ch$-module $\Hom_P(\EE,V^{U_P})$. Let $\g\in G(K)_{cvr}$ be such that 
$\g\in Q=P_H$. Let $T'=T'_\g$. Let $\bg$ be the image of $\g$ in $Q/U_Q$. Note that $\bg$ is an 
$F_q$-rational regular semisimple element of $\uQ/U_{\uQ}$, see 2.4(b). Let $\boc$ be the conjugacy class in
$W_H$ defined as follows: if $\b'$ is a Borel subgroup of $\uQ/U_{\uQ}$ containing $\bg$ and $\b''$ is the 
image of $\b'$ under the Frobenius map then the relative position of $\b',\b''$  belongs to $\boc$. Note 
that $\boc\sub C_{T'}$ where $C_{T'}$ is the conjugacy class in $W'$ associated to $T'$ in 2.5(a). For any 
$w\in\boc$ let $R_w$ be the virtual representation of the finite reductive group $Q/U_Q$ associated to $w$ 
as in \cite{\DL}. Since $R_w$ is independent of $w$ as long as $w\in\boc$, we denote it by $R_\boc^H$.
According to \cite{\DL, 7.9} for any unipotent representation $E$ of $Q/U_Q$ we have 
$$\tr(\bg,E)=(E:R^H_\boc)$$
where $(:)$ denotes multiplicity. Now 3.6 implies
$$\tr(\g,V^{U_Q})=\sum_{i\in[1,r]}\tr(\bg,[D_i])=\sum_{i\in[1,r]}([D_i]:R_\boc^H)\in\ZZ.\tag a$$

\subhead 3.8\endsubhead
For $V,P$ as above and for $Q=P_H$ with $H\sneq\II$ (but without assuming that $P\sub Q$) we show:

(a) {\it If $V^{U_Q}\ne0$ then $J\sub\o(H)$ for some $\o\in\Om$.}
\nl
(Note that $\Om$ (see 2.1) acts by conjugation on $\{s_i;i\in\II\}$ and this induces an action of $\Om$ on 
$\II$.) The proof is almost a repetition of that of 3.5(a).

We can find a parabolic subgroup $\p'$ of $Q/U_Q$ and an irreducible cuspidal representation $\EE_0$ of 
$\p'/u_{\p'}$ such that the $\p'/u_{\p'}$-module $E_1^{u_{\p'}}$ contains $\EE_0$. Then the 
$\p'/u_{\p'}$-module $(V^{U_Q})^{u_{\p'}}$ contains $\EE_0$. Let $P'$ be the inverse image of $\p'$ under 
$Q@>>>Q/U_Q$; then $P'\in\cp$ is contained in $Q$ and the $P'/U_{P'}$-module $(V^{U_Q})^{U_{P'}}$ contains 
$\EE_0$. (We have $P'/U_{P'}=\p'/u_{\p'}$.) Replacing $P'$ by a $Q$-conjugate we can assume that
$P'=P_{J_1}$ for some $J_1\sub H$. We have $U_Q\sub U_{P'}$ hence $V^{U_{P'}}\sub V^{U_Q}$ and  
$(V^{U_Q})^{U_{P'}}=V^{U_{P'}}$. Thus, the $P'/U_{P'}$-module $V^{U_{P'}}$ contains $\EE_0$. By the 
uniqueness of cuspidal support of $V$ we see that $\EE_0$ is a unipotent (cuspidal) representation and then 
using \cite{\CLA, 1.6(b)} we see that we can find $g'\in G(K)$ such that $\Ad(g')$ carries $P$ to $P'$ and 
$\EE$ to a representation isomorphic to $\EE_0$. In particular we have $\Ad(g')(P)\sub Q$, proving (a).

\head 4. The character formula\endhead
\subhead 4.1\endsubhead
Let $T$ be as in 1.2 and let $B$ be as in 3.1. Let $J\sneq\II$, $P=P_J$, $\EE$ be as in 3.2 and let $(V,\s)$
be as in 3.3. Let $\ph_V$ be as in 1.1. 

Let $\g\in G(K)_{cvr}$. We want to compute $\ph_V(\g)$. As in 2.7(b) we associate to $\g$ a subset
$H\sneq\II$ so that (after replacing if necessary $\g$ by a $G(K)_{der}$-conjugate), $Q:=P_H$ contains $\g$ 
and (setting $T'=T'_\g$) we have $T'_{*s}\sub T_*$ (notation of 2.7). Let $M$ be the centralizer of 
$T'_{*s}$, see 2.7; this is the Levi subgroup of a parabolic subgroup of $G$. Let $T'_0$ (resp. $T'_1$) be 
the maximal compact subgroup (resp. the maximal pro-p subgroup) of $T'(K)$. Let $Q=Q_0,Q_1,Q_2,\cdots$ be 
the strictly decreasing Moy-Prasad filtration \cite{\MP} of $Q$. Then $Q_1=U_Q$. By \cite{\AK}, \cite{\MS},

(a) {\it $\ph_V$ is constant on $\cup_{x\in G(K)}x\g T'_1x\i$.}

\proclaim{Proposition 4.2} $\g Q_1\sub\cup_{x\in G(K)}x\g T'_1x\i$.
\endproclaim
For $n\ge0$ we set $\cm_n:=M\cap Q_n$, $T'_n=T'\cap Q_n$. (This agrees with the earlier definition of 
$T'_0,T'_1$.) Note that $Q_0/Q_1\simeq\cm_0/\cm_1$ and $Q_i/Q_{i+1}$ is abelian for $i\ge1$. Moreover, 
$\cm_0/\cm_1$ (resp. $T'_0/T'_1$) is the group of the $F_q$-points of a reductive group over $F_q$ whose 
Lie algebra is denoted by $\fm$ (resp. $\ft'$). 

We need the following two lemmas:

\proclaim{Lemma 4.3} Let $n\ge1$ and $\g'\in\g T'_1$. Then, for any $z_n\in\cm_n$, there exist 
$g_n\in\cm_n$, $\g_n\in T'_n$ such that $g_n\g'z_ng_n^{-1}\in\g'\g_n\cm_{n+1}$. 
\endproclaim
If $\cm_n=\cm_{n+1}$, one can take $z_n=1=g_n$. If $\cm_n\ne\cm_{n+1}$, there is a natural isomorphism 
$\iota_n:\cm_n/\cm_{n+1}\rightarrow\ft'(F_q)\oplus(\fm/\ft')(F_q)$.  Note that since $\g'\in G(K)_{cvr}$, 
$d:=\Ad(\g'{}\i)-1$ induces an isomorphism $(\fm/\ft')(F_q)\rightarrow(\fm/\ft')(F_q)$. 
Let $z_n\in\cm_n$ and $\iota_n(z_n)=t_n+m_n$ with $t_n\in\ft'(F_q)$ and $m_n\in(\fm/\ft')(F_q)$. 
Choose $\g_n\in T'_n$ and $g_n\in\cm_n$ so that $\iota_n(\g_n)=t_n$ and $\iota_n(d(g_n))=-m_n$. Then 
$\g'{}\i g_n\g'z_ng_n\i\in\g_n\cm_{n+1}$. Setting $z_{n+1}=\g_n^{-1}\g'{}\i g_n\g'z_ng_n^{-1}$, we have 
$g_n\g'z_ng_n^{-1}=\g'\g_nz_{n+1}\in\g'\g_n\cm_{n+1}$.

\proclaim{Lemma 4.4} Let $n\ge1$. For any $\g'\in\g T'_1$ and $z\in Q_n$, there exist $g\in Q_n$, 
$\g_n\in T'_n$ and $z'\in Q_{n+1}$ such that $\Ad(g)(\g'z)=\g'\g_nz'$.
\endproclaim
Let $Z=\{\a\in R;U_\a\cap Q_n\supsetneq U_\a\cap Q_{n+1}\}$. If $Z=\emptyset$, then $Q_n=\cm_nQ_{n+1}$. 
Hence $z=z_nz'$ for some $z_n\in\cm_n$, $z'\in Q_{n+1}$ and one can take $g=g_n$ and $\g_n$ as in Lemma 4.3. 
If $Z\neq\emptyset$, then $\cm_n=\cm_{n+1}$ and $T'_n=T'_{n+1}$. Moreover, since $\g'\in G(K)_{cvr}$, 
$d:=\Ad(\g'{}\i)-1$ induces an isomorphism $Q_n/Q_{n+1}\rightarrow Q_n/Q_{n+1}$. 
Let $g\in Q_n$ be such that $d(g)=z\i$. Then, $\g'{}\i g\g'zg\i\in Q_{n+1}$. 
Setting $z'=\g'{}\i g\g'zg\i$ and $\g_n=1$, we have $g\g'zg\i=\g'\g_nz'\in\g'\g_nQ_{n+1}$.

\mpb

Continuing with the proof of Proposition 4.2, let $\g_0=\g$ and $z_1\in Q_1$. We will construct inductively 
sequences $g_1,g_2,\cdots$, $t_1,t_2,\cdots$ and $z_1,z_2,\cdots$ such that $g_i\in Q_i$, $\g_i\in T'_i$, 
$z_i\in Q_i$ and
$$\Ad(g_k\cdots g_2g_1)(\g_0z_1)=\Ad(g_k)(\g_0\g_1\cdots\g_{k-1}z_k)=(\g_0\g_1\cdots \g_k)z_{k+1}.$$
Applying Lemma 4.4 with $n=1$, $\g'=\g_0$ and $z=z_1$, we find $\g_1\in T'_1$ and $z_2\in Q_2$ such that 
$g_1\g_0z_1g_1^{-1}=\g_0\g_1z_2$ with $\g_1\in T'_1$ and $z_2\in Q_2$. 
Suppose that we have found $g_i\in Q_i$, $z_{i+1}\in Q_{i+1}$ and $\g_i\in T'_i$ for $i=1,\cdots,k$ where 
$k\ge1$. Applying Lemma 4.4 with $n=k+1$, $\g'=\g_0\g_1\cdots\g_k$ and $z=z_{k+1}$, we find 
$g_{k+1}\in Q_{k+1}$, $\g_{k+1}\in T'_{k+1}$ and $z_{k+2}\in Q_{k+2}$ such that 
$$g_{k+1}\g_0\g_1\cdots\g_kz_{k+1}g_{k+1}^{-1}=\Ad(g_{k+1}\cdots g_2g_1)(\g_0z_1)
=\g_0\g_1\g_2\cdots\g_{k+1}z_{k+2}.$$
Taking $g\in Q_1$ to be the limit of $g_k\cdots g_2g_1$ as $k\rightarrow\infty$, we have 
$\Ad(g)(\g z_1)\in\g T'_1$. This completes the proof of Proposition 4.2.

Note that $\g$ acts on $V^{Q_1}=V^{U_Q}$ since $\g$ normalizes $Q_1$.

\proclaim{Proposition 4.5} $\ph_V(\g)=\tr(\g,V^{U_Q})$.
\endproclaim
If $f:G(K)@>>>\CC$ is any locally constant function with compact support then there is a well defined linear
map $\s_f:V@>>>V$ such that for any $x\in V$ we have $\s_f(x)=\int_Gf(g)\s(g)(x)dg$. ($dg$ is the Haar 
measure on $G(K)$ for which $vol(B)=1$). This linear map has finite rank hence it has a well defined trace 
$\tr(\s_f)\in\CC$.

Note that $Q_1\g Q_1=\g Q_1$. Let $f_\g$ be the characteristic function of $\g Q_1$. Then, by 4.1(a) and 
Proposition 4.2 we have
$$\align&\tr({\s_{f_\g}})=\int_G f_\g(g)\ph_V(g)dg\\&
=\int_{Q_1\g Q_1}\ph_V(\g)\,dg=vol(Q_1\g Q_1)\ph_V(\g)=vol(Q_1)\ph_V(\g).\endalign$$
Moreover, $\s_{f_\g}$ maps $V$ onto $V^{Q_1}$, and on $V^{Q_1}$, it is equal to $\s(\g)$ times $vol(Q_1)$. 
Hence $\tr({\s_{f_\g}})=\tr({\s_{f_\g}};V^{Q_1})=vol(Q_1)\tr(\s(\g)|V^{Q_1})$ and the proposition follows.

\mpb

Note that the proposition above is in fact valid (with the same proof) for any irreducible admissible
(not necessarily unipotent) representation $V$ of depth $0$ of $G(K)$. We now state the main result of this 
paper. 

\proclaim{Theorem 4.6} We preserve the setup of 4.1. 

(a) If $J\not\sub\o(H)$ for any $\o\in\Om$ then $\ph_V(\g)=0$.

(b) Assume that $J\sub H$. Then 
$$\ph_V(\g)=\sum_{i\in[1,r]}([D_i]:R_\boc^H)\in\ZZ$$
(notation of 3.7).
\endproclaim
We prove (a). Assume that $\ph_V(\g)\ne0$. By 4.5 we then have $V^{U_Q}\ne0$ hence by 3.8(a), 
$P_J\sub P_{\o(H)}$ for some $\o\in\Om$. It follows that $J\sub\o(H)$ and (a) is proved. Now (b) is obtained
by combining 4.5 with 3.7(a).

\subhead 4.7\endsubhead
The case where $J\sub\o(H)$ for some $\o\in\Om$ can be reduced to the case 4.6(b) as follows. Let $\tio$ be 
an element in the normalizer of $B$ in $G(K)$ that represents $\o$; we can assume that $\tio$ normalizes 
$T$. Let $\g'=\tio\i\g\tio,\ti Q=\tio\i Q\tio\i=P_{\o\i(H)}$. Now 4.6(b) is applicable to $\g',\ti Q$ 
instead of $\g,Q$ and it yields a formula for $\ph_V(\g')$. Since $\g,\g'$ are conjugate in $G(K)$ we have 
$\ph_V(\g)=\ph_V(\g')$ so that we get a formula for $\ph_V(\g)$.

\subhead 4.8\endsubhead
Let $C$ be a conjugacy class of elements of finite order in $W'$ and let $\g',\g''$ be two elements in 
$G(K)_{cvr}$ such that (denoting by $T'$, $T''$ the maximal tori in $G$ that contain $\g',\g''$) we have
$C_{T'}=C_{T''}=C$. We show:

(a) $\ph_V(\g')=\ph_V(\g'')$.
\nl
(It follows that $\ph_V(\g)$ is independent of $\g$ as long as $\g\in G(K)_{cvr}$ gives rise to a fixed
conjugacy class in $W'$.)

We can assume that $\g'=\g$ and $Q=P_H$ are as in 4.1. We can also assume that $T''=T'$ so that $\g''\in T'$.
Using 2.8 we see that $Q$ is a minimal parahoric subgroup of $G(K)$ containing $\g''$. Using 4.6 and 4.7 we 
now see that (a) holds.

\head 5. Examples and comments\endhead
\subhead 5.1\endsubhead
In this and the next two subsections we assume that $G$ is adjoint and that 4.1(a) holds. We denote by $G^*$
a semisimple group over $\CC$ of type dual to that of $G$. The unipotent representations of $G(K)$ were 
classified in \cite{\CLA} (extending the classification in \cite{\KL} of Iwahori-spherical representations).
Namely in \cite{\CLA, 5.21} a bijection is established between the set $\cu$ of unipotent representations of
$G(K)$ (up to isomorphism) and the set $\fS$ consisting of all triples $(s,u,\r)$ (modulo the natural 
action of $G^*$) where $s$ is a semisimple element of $G^*$, $u$ is a unipotent element of $G^*$ such that 
$sus\i=u^q$ and $\r$ is an irreducible representation (up to isomorphism) of the group $\bZ(s,u)$ of 
components of the centralizer $Z(s,u)$ of $s$ and $u$ in $G^*$ such that the centre of $G^*$ acts trivially 
on $\r$); let $V_{s,u,\r}$ be the unipotent representation corresponding to $(s,u,\r)$.

In the case where $V=V_{s,u,\r}$ in Theorem 4.6 is tempered, the definition of $V$ in terms of equivariant 
homology given in \cite{\CLA} shows that the restriction of the $\ch$-module $\Hom_P(\EE,V^{U_P})$ (see 3.3) 
to the subalgebra $\ch_H$ is explicitly computable (in terms of generalized Green functions); moreover, the 
multiplicities $([D_i]:R_\boc^H)$ in 4.6(b) are also explicitly known from \cite{\ORA, 4.23} (in terms of 
the nonabelian Fourier transform of \cite{\ORA, 4.14}). We see that in this case the character values 
$\ph_V(\g)$ in 4.6 are explicitly computable integers.

\subhead 5.2\endsubhead
In this and next two subsections we assume that $G$ is of type $G_2,F_4$ or $E_8$.

For any finite group $\G$ let $M(\G)$ be the set of all pairs $(y,r)$ where $y$ is an element of 
$\G$ defined up to conjugacy and $r$ is an irreducible representation of the centralizer $Z_\G(y)$ of 
$y$ in $\G$, up to isomorphism.

We have a partition $\fS=\sqc_{\fc}\fS_\fc$ where $\fc$ runs over the unipotent conjugacy classes in $G^*$ 
and $\fS_\fc$ is the set of (equivalence classes of) triples $(s,u,\r)$ in $\fS$ with $u\in\fc$.
Let $\cu_\fc$ be the subset of $\cu$ corresponsing to $\fS_\fc$ under the bijection $\cu\lra\fS$ in 5.1.

Now assume that $\fc$ is a distinguished unipotent class in $G^*$ and let $u\in\fc$. Let $\G$ (resp. $\G'$) 
be the group of connected components of $\ti\G=\{g\in G^*;gug\i=g\}$ (resp. of 
$\ti\G'=\{(g,\l)\in G^*\T\CC^*;gug\i=u^\l\}$). The obvious imbedding $\ti\G@>>>\ti\G'$ induces an 
isomorphism $\G@>\si>>\G'$. 

If $(s,u,\r)$ represents an element of $\fS_\fc$ then the connected
component of $(s,q)$ in $\ti\G'$ can be identified with a connected component of $\ti\G$ hence with an
element $y\in\G$; moreover $\r$ can be viewed as an irreducible representation $r$ of $Z_\G(y)$. Thus 
$\fS_\fc$ is identified with $M(\G)$. We shall write $V_{\fc;y,r}$ instead of $V_{s,u,\r}$ when $(y,r)$ 
corresponds as above to $(s,u,\r)\in\fS_\fc$.
Let $M'(\G)$ be the set of all $(y,r)\in M(\G)$ such that $V=V_{\fc;y,r}$ has a nonzero space of
$B$-invariants $V^B$. For such $(y,r)$, $V^B$ is a square integrable irreducible representation of the
ordinary affine Hecke algebra $\ch$. By the $K$-theoretic construction of this representation \cite{\KL},
$V^B$ is the specialization at $v=\sqrt{q}$ of a representation of the affine Hecke algebra with parameter 
$v$; this representation can be specialized at $v=1$, yielding a (finite dimensional)
representation $A_{y,r}$ of $W$ (which in particular is a representation of $W'$).
For $(y,r),(y',r')$ in $M(\G)$ let $\{(y,r),(y',r')\}\in\CC$ be the 
$((y,r),(y',r'))$-entry of the nonabelian Fourier transform matrix on $M(\G)$ (see \cite{\ORA, 4.14}).

It seems likely that in our case with $V\in\cu_\fc$, the character formula in 4.6 can be rewritten in the 
following form (we write $V=V_{\fc;y,r}$ where $(y,r)\in M(\G)$):
$$\ph_V=\sum_{(y',r')\in M'(\G)}\{(y,r),(y',r')\}\t_{A_{y',r'}}\tag a$$
as functions $G(K)_{cvr}@>>>\ZZ$, where $\t_{A_{y',r'}}$ is as in 2.3.

For example, (a) holds when $G$ is of type $E_8$ and $\fc$ is the regular, or subregular or 
subsubregular unipotent class in $G^*$. (In these cases we have $\G=\{1\}$.)

\subhead 5.3\endsubhead
In this subsection we assume that $G$ is of type $G_2$ and $\fc$ is the subregular unipotent class in $G^*$.
In this case we have $\G=S_3$ and $\cu_\fc$ consists of eight square integrable irreducible representations 
of $G(K)$ considered in \cite{\SQ, 1.7}. Four representations in $\cu_\fc$ (denoted by 
$V,V',V'',V'''$) are Iwahori-spherical and the four corresponding
irreducible representations of the affine Hecke algebra of type $G_2$ are carried by the $W$-graph denoted in
\cite{\SQ, 3.13, Type $\ti G_2$} by $\cg,\cg',\cg'',\cg'''$ (they have dimension $3,3,2,1$ and $\cg$ 
gives rise to the reflection representation of the affine Hecke algebra).
The other four representations in $\cu_\fc$ are supercuspidal and are denoted by $S,S',S'',S'''$.
The $W$-graphs $\cg,\cg',\cg'',\cg'''$ also carry (for $q=1$) representations of $W$ of dimension $3,3,2,1$
denoted by $A,A',A'',A'''$. From Theorem 4.6 we can deduce:
$$\ph_V=(1/6)\t_{A}+(1/2)\t_{A'}+(1/3)\t_{A''}+(1/3)\t_{A''''},$$
$$\ph_{V'}=(1/2)\t_A+(1/2)\t_{A'},$$
$$\ph_{V''}=(1/3)\t_A +(2/3)\t_{A''}-(1/3)\t_{A'''},$$
$$\ph_{V'''}=(1/3)\t_A -(1/3)\t_{A''}+(2/3)\t_{A'''},$$
$$\ph_{S}=(1/6)\t_A-(1/2)\t_{A'}+(1/3)\t_{A''}+(1/3)\t_{A'''},$$
$$\ph_{S'}=(1/2)\t_A-(1/2)\t_{A'},$$
$$\ph_{S''}=\ph_{S'''}=(1/3)\t_A -(1/3)\t_{A''}-(1/3)\t_{A'''},$$
as functions $G(K)_{cvr}@>>>\ZZ$.
We see that in our case the equality 5.2(a) holds.

\subhead 5.4\endsubhead
The functions $\t_A$ for $A=A_{y,r}$ ($(y,r)\in M'(\G)$ in the setup of 5.2)
can be regarded as $p$-adic analogues of the (uniform) almost characters (see \cite{\ORA}) for the 
corresponding reductive group over $F_q$ (restricted to $F_q$-rational regular semisimple elements).
It would be interesting to see whether these functions have natural extensions to
$G_{rs}\cap G(K)$ (in the same way as the uniform almost characters of reductive groups over $F_q$ are 
defined on the whole group not just on regular semisimple elements) and then to see whether these
extensions play a role in computing $\ph_V$ on $G_{rs}\cap G(K)$.

\subhead 5.5\endsubhead
Let us now drop the assumption in 1.1 that $G$ is split over $K$ (but we still assume that $G$ is split over 
$\uK$). Then the unipotent representations of $G(K)$ are well defined, see \cite{\CLA}, and $G(K)_{cvr}$ is 
defined as in 1.2. We will show elsewhere that the character formula 4.6 for 
unipotent representations of $G(K)$ extends to this more general case with 
essentially the same proof.

\

{\it Note added 3/06/2013.} We thank the referee for pointing out that Lemmas 4.3, 4.4 and Propositions 
4.2, 4.5 are special cases of Proposition 7.4, Lemma 7.5 and Theorem 8.1 of \cite{\MS}. In an earlier 
version of this paper, the results of this paper were proved under some restrictions on characteristic, so 
that we could refer to \cite{\AK}. We thank the referee for pointing out the results of \cite{\MS} which 
allow us to remove the restrictions on characteristic.

An extension of the type mentioned in 5.4 has meanwhile been found in
[G. Lusztig, Unipotent almost characters of simple $p$-adic groups,
arxiv:1212.6540].

\enddocument